\documentclass[12pt]{amsart}
\usepackage{amssymb,latexsym}
\usepackage{pdfsync}

\newdimen\AAdi%
\newbox\AAbo%
\def\AAk#1#2{\s_etbox\AAbo=\hbox{#2}\AAdi=\wd\AAbo\kern#1\AAdi{}}%
\def\AAr#1#2#3{\s_etbox\AAbo=\hbox{#2}\AAdi=\ht\AAbo\raise#1\AAdi\hbox{#3}}%
\font\tenmsb=msbm10 at 12pt
\font\sevenmsb=msbm7 at 8pt
\font\fivemsb=msbm5 at 6pt
\newfam\msbfam
\textfont\msbfam=\tenmsb
\scriptfont\msbfam=\sevenmsb
\scriptscriptfont\msbfam=\fivemsb
\def\Bbb#1{{\tenmsb\fam\msbfam#1}}

\newcommand{\beq}{\begin{equation}}
\newcommand{\eeq}{\end{equation}}
\newcommand{\beqr}{\begin{eqnarray}}
\newcommand{\eeqr}{\end{eqnarray}}
\newcommand{\ba}{\begin{array}}
\newcommand{\ea}{\end{array}}

\begin{document}

\newtheorem{thm}{Theorem}
\newtheorem{lem}{Lemma}
\newtheorem{cor}{Corollary}
\newtheorem{rem}{Remark}
\newtheorem{pro}{Proposition}
\newtheorem{defi}{Definition}
\newtheorem{eg}{Example}
\newtheorem*{claim}{Claim}
\newtheorem{conj}[thm]{Conjecture}
\newcommand{\noi}{\noindent}
\newcommand{\dis}{\displaystyle}
\newcommand{\mint}{-\!\!\!\!\!\!\int}
\numberwithin{equation}{section}

\def \bx{\hspace{2.5mm}\rule{2.5mm}{2.5mm}}
\def \vs{\vspace*{0.2cm}}
\def\hs{\hspace*{0.6cm}}
\def \ds{\displaystyle}
\def \p{\partial}
\def \O{\Omega}
\def \o{\omega}
\def \b{\beta}
\def \m{\mu}
\def \l{\lambda}
\def\L{\Lambda}
\def \ul{u_\lambda}
\def \D{\Delta}
\def \d{\delta}
\def \k{\kappa}
\def \s{\sigma}
\def \e{\varepsilon}
\def \a{\alpha}
\def \tf{\tilde{f}}
\def\cqfd{%
\mbox{ }%
\nolinebreak%
\hfill%
\rule{2mm} {2mm}%
\medbreak%
\par%
}
\def \pr {\noindent {\it Proof.} }
\def \rmk {\noindent {\it Remark} }
\def \esp {\hspace{4mm}}
\def \dsp {\hspace{2mm}}
\def \ssp {\hspace{1mm}}

\def\la{\langle}\def\ra{\rangle}

\def \u{u_+^{p^*}}
\def \ui{(u_+)^{p^*+1}}
\def \ul{(u^k)_+^{p^*}}
\def \energy{\int_{\R^n}\u }
\def \sk{\s_k}
\def \mo{\mu_k}
\def\cal{\mathcal}
\def \I{{\cal I}}
\def \J{{\cal J}}
\def \K{{\cal K}}
\def \OM{\overline{M}}

\def\n{\nabla}

\def\fk{{{\cal F}}_k}
\def\M1{{{\cal M}}_1}
\def\Fk{{\cal F}_k}
\def\Fl{{\cal F}_l}
\def\FF{\cal F}
\def\Gk{{\Gamma_k^+}}
\def\n{\nabla}
\def\uuu{{\n ^2 u+du\otimes du-\frac {|\n u|^2} 2 g_0+S_{g_0}}}
\def\uuug{{\n ^2 u+du\otimes du-\frac {|\n u|^2} 2 g+S_{g}}}
\def\sku{\sk\left(\uuu\right)}
\def\qed{\cqfd}
\def\vvv{{\frac{\n ^2 v} v -\frac {|\n v|^2} {2v^2} g_0+S_{g_0}}}
\def\vvs{{\frac{\n ^2 \tilde v} {\tilde v}
 -\frac {|\n \tilde v|^2} {2\tilde v^2} g_{S^n}+S_{g_{S^n}}}}
\def\skv{\sk\left(\vvv\right)}
\def\tr{\hbox{tr}}
\def\pO{\partial \Omega}
\def\dist{\hbox{dist}}
\def\RR{\Bbb R}\def\R{\Bbb R}
\def\C{\Bbb C}
\def\B{\Bbb B}
\def\N{\Bbb N}
\def\Q{\Bbb Q}
\def\Z{\Bbb Z}
\def\PP{\Bbb P}
\def\EE{\Bbb E}
\def\F{\Bbb F}
\def\G{\Bbb G}
\def\H{\Bbb H}
\def\SS{\Bbb S}\def\S{\Bbb S}

\def\div{\hbox{div}\,}

\def\lcf{{locally conformally flat} }

\def\circledwedge{\setbox0=\hbox{$\bigcirc$}\relax \mathbin {\hbox
to0pt{\raise.5pt\hbox to\wd0{\hfil $\wedge$\hfil}\hss}\box0 }}

\def\sss{\frac{\s_2}{\s_1}}

\date{\today}
\title[ A Bernstein type result of translating solitons ]{ A Bernstein type result of translating solitons }

\author{}

\author[Qiu]{Hongbing Qiu}
\address{School of Mathematics and Statistics\\ Wuhan University\\Wuhan 430072,
China,
and Hubei Key Laboratory of Computational Science \\ Wuhan University\\Wuhan 430072,
China
 }
 \email{hbqiu@whu.edu.cn}

 \thanks{
 This work is partially supported by NSFC (No. 11771339) and Hubei Provincial Natural Science Foundation of China (No. 2021CFB400).  }

\begin{abstract}
 
By using a new test function and the gradient estimate technique, we obtain a better Bernstein type result of translating solitons.

\vskip12pt



\end{abstract}
\maketitle
\section{Introduction}

Let $X: M^{n} \rightarrow \mathbb{R}^{m+n}$ be an isometric immersion from an $n$-dimensional oriented Riemannian manifold $M$ to the Euclidean space $\mathbb{R}^{m+n}$.  
The mean curvature flow (MCF) in Euclidean space is a one-parameter family of immersions $X_t= X(\cdot, t): M^m \rightarrow \mathbb{R}^{m+n}$ with the corresponding image $M_t=X_t(M)$ such that
\begin{equation}\label{eqn-MCF}
\begin{cases}\aligned
\frac{\partial}{\partial t}X(x,t)=& H(x,t), x\in M,\\
X(x,0)=&X(x),
\endaligned
\end{cases}
\end{equation}
is satisfied, where $H(x, t)$ is the mean curvature vector of $M_t$ at $X(x, t)$ in $\mathbb{R}^{m+n}$.

We call
$M^n$ a translating soliton in $\mathbb{R}^{m+n}$ if it
 satisfies
\begin{equation}\label{eqn-T111}
H= V_{0}^N,
\end{equation}
where  $V_0$ is a fixed vector in $\mathbb{R}^{m+n}$ with unit length and $V_{0}^N$ denotes the orthogonal projection of
$V_0$ onto the normal bundle of $M^n$. The translating solitons give rise to eternal solutions $X_t=X+tV_0$ to (\ref{eqn-MCF}). They are not only  special solutions to the mean curvature flow equations, but also they often occur as the Type ${\rm II}$ singularity of mean curvature flows (see \cite{AV95, AV97, Ham95, HuiSin99, HuiSin2-99, Whi00, Whi03}). And the geometry of the translating soliton has been paid much attention during the past two decades, see the references (not exhaustive): \cite{CheQiu16, CSS, Hal, Jian, Kun, MSS, Ngu09, Ngu13, NT, Qiu, Wang11, Xin}, etc.

Note that the minimal submanifold is a static solution of the MCF, we firstly recall the related theory of minimal submanifolds. 
The Bernstein problem has been a central problem in the study of minimal submanifolds. The classical Bernstein theorem states that the entire minimal graph in the Euclidean space $\mathbb{R}^3$ is a plane \cite{Bern}.
Many efforts had been made to generalize the Bernstein theorem to higher dimensions. Eventually, Simons \cite{Sim} proved that entire minimal graphic hypersurfaces in $\mathbb{R}^{n+1}$ must be hyperplanes for $n\leq 7$ (see De Giorgi \cite{DG} for $n=3$ and Almgren \cite{Alm} for $n=4$), while Bombieri-De Giorgi-Giusti \cite{BDG} constructed a counterexample for $n\geq 8$. However, Moser \cite{Mos} had earlier showed that, under an additional assumption that the gradient of the graph function is uniformly bounded, the entire minimal graphic hypersurface has to be planar in arbitrary dimension (see also \cite{HJW}). Later, Ecker-Huisken \cite{EH} improved Moser’s theorem by using the curvature estimate technique. There are plenty of works on the study of the Bernstein type problems on minimal hypersurfaces in Euclidean spaces (see  \cite{CP, FS, SSY, Sol, JosXinYan12, Din20, Oss, Xav, Fuj}).

The Bernstein problem in higher codimensions becomes more complicated, for example, the Moser's theorem can be generalized to higher codimensional complete minimal graphs for dimension 2 and 3 (see \cite{CO, Bar, Fis}), but the counterexample of a nontrivial minimal graph with bounded slope, constructed by Lawson-Osserman \cite{LO}, 
 sets a limit for how far one can go. They also raised a question for finding the ``best" constant possible in the same paper. Afterward, the Bernstein type theorem can be achieved that any minimal graph $M^n$ in $\mathbb{R}^{m+n} (n\geq 3, m\geq 2)$ must be an affine $n$-plane provided the slope $\leq 3$ (see \cite{HJW, JX, XY, JXY, JXY18}). Recently, Assimos-Jost \cite{AJ} extended Moser's theorem to codimension 2 by using Sampson's maximum principle.
 
 Since the translating soliton can be viewed as a generalization of minimal submanifolds, it is natural to study the Bernstein type problem of translating solitons. For the codimension one case, Bao-Shi \cite{BS} proved a translating soliton version of Moser's theorem, namely, if the image of the translating soliton $M^n$ under the Gauss map is contained in a regular ball, then such a complete translating soliton in $\mathbb{R}^{n+1}$ has to be a hyperplane. For higher codimensions, Kunikawa \cite{Kun} generalized the result of \cite{BS} to the flat normal bundle case, and obtained that for any complete translating soliton $M^n$ with flat normal bundle in $\mathbb{R}^{m+n}$,  if the $w$-function is positive and it satisfies the growth condition $w^{-1} = o(R^{\frac{1}{2}})$, then $M^n$ must be an affine subspace. 
 In general, without the condition on the flat normal bundle, Xin \cite{Xin} partially solved a counterpart of the Lawson-Osserman problem for translating solitons, that is, if the $v$-function satisfies $v\leq v_1 < v_0:= \frac{2\cdot 3^{\frac{2}{3}}}{1+3^{\frac{2}{3}}}$, then any complete translating soliton $M^n$ in $\mathbb{R}^{m+n}_n (m\geq 2)$ has to be affine linear. It is natural to ask that whether we can find a larger constant $v_0$, such that the corresponding Bernstein type result still holds.

In this note, by adopting a new test function which is different from the one in \cite{Xin} and using the gradient estimate technique, we have the following result.

\begin{thm}\label{thm1}

Let $M^n$ be a complete $n$-dimensional translating soliton in $\mathbb{R}^{m+n}$ with codimension $m \geq 2$ and the positive $w$-function. Put $v_0:= \frac{2\cdot 4^{\frac{2}{3}}}{1+ 4^{\frac{2}{3}}}$. If for any constant $v_1 < v_0$, the $v$-function satisfies 
\[
v\leq v_1 < v_0,
\]
then $M^n$ is affine linear.

\end{thm}

\begin{rem}

Xin \cite{Xin} showed that the same conclusion holds under the condition that $v\leq v_1<\frac{2\cdot 3^{\frac{2}{3}}}{1+ 3^{\frac{2}{3}}}$. Clearly, $\frac{2\cdot 3^{\frac{2}{3}}}{1+ 3^{\frac{2}{3}}}<\frac{2\cdot 4^{\frac{2}{3}}}{1+ 4^{\frac{2}{3}}}$. Thus the condition in Theorem \ref{thm1} is weaker than the one in \cite{Xin}.

\end{rem}

\begin{rem}

By Proposition 6.1 in \cite{XY} (see also \cite{XY09}), if a $v$-function is bounded from above by a positive constant, then any complete translating soliton $M^n$ is an entire graph. Let $u^\a: \mathbb{R}^n \to \mathbb{R}^m$ be the graph functions. Then the induced metric $(g_{ij})$ on $M^n$ is $g_{ij} = \delta_{ij} + \sum_\a u^\a_{i}u^\a_{j}$ and the $v$-function is just $\left( \det\left(  \delta_{ij} + \sum_\a u^\a_{i}u^\a_{j} \right) \right)^{\frac{1}{2}}$. Theorem \ref{thm1} claims that any entire graphic translating soliton $M^n$ in $\mathbb{R}^{m+n} (m\geq 2)$ has to be affine linear provided
\[
\left( \det\left(  \delta_{ij} + \sum_\a u^\a_{i}u^\a_{j} \right) \right)^{\frac{1}{2}} \leq v_1 < v_0:= \frac{2\cdot 4^{\frac{2}{3}}}{1+ 4^{\frac{2}{3}}}.
\]

\end{rem}

\vskip24pt

\section{Preliminaries}

Let $G_{n, m}$ be the Grassmann manifold consisting of the oriented linear $n$-subspaces in $\mathbb{R}^{m+n}$. The canonical Riemannian structure on $G_{n, m}$ makes it a natural generalization of teh Euclidean sphere. $G_{n, m}=SO(m+n)/SO(n)\times SO(m)$ is an irreducible symmetric space of compact type.

For every $P \in G_{n, m}$, we choose an oriented basis $\{ u_1, \cdot \cdot \cdot, u_n \}$ of $P$, and let
\[
\psi(P):= u_1 \wedge \cdot\cdot\cdot \wedge u_n \in \Lambda^n (\mathbb{R}^{m+n}).
\]
A different basis for $P$ shall give a different exterior product, but the two products differ only by a positive scalar; $\psi(P)$ is called the Pl${\ddot{\rm u}}$cker coordinate of $P$, which is a homogeneous coordinate.

 Via the Pl${\ddot{\rm u}}$cker embedding, $G_{n, m}$ can be viewed as a submanifold of some Euclidean space. The restriction of the Euclidean inner product on $G_{n, m}$ is denoted by $w: G_{n, m} \times G_{n, m} \to \mathbb{R}$
\[
w(P, Q) = \frac{\la \psi(P), \psi(Q) \ra}{\la \psi(P), \psi(P) \ra^{\frac{1}{2}}\la \psi(Q), \psi(Q) \ra^{\frac{1}{2}}}.
\] 
If $\{ e_1, \cdot\cdot\cdot, e_n \}$ is an oriented orthonormal basis of $P$ and $\{ f_1, \cdot\cdot \cdot, f_n\}$ is an oriented orthonormal basis of $Q$, then 
\[
w(P, Q) = \la e_1 \wedge\cdot\cdot\cdot \wedge e_n, f_1\wedge \cdot\cdot\cdot \wedge f_n \ra = \det W,
\]
where $W=(\la e_i, f_j \ra)$. It is well known that 
\[
W^TW=O^T\Lambda O
\]
with $O$ an orthogonal matrix and $\Lambda = {\rm diag}(\mu_1^{2}, \cdot\cdot\cdot, \mu_n^{2})$. Here each $0\leq \mu_i^{2}\leq 1$. Putting $p:= \min\{ m, n \}$, then at most $p$ elements in $\{ \mu_1^{2}, \cdot\cdot\cdot, \mu_n^{2} \}$ are not equal to 1. Without loss of generality, we can assume $\mu_i^{2} = 1$ whenever $i>p$. We also note that the $\mu_i^{2}$ can be expressed as 
\[
\mu_i^{2} = \frac{1}{1+\lambda_i^{2}}
\]
with $\lambda_i \in [0, +\infty)$.

The Jordan angles between $P$ and $Q$ are defined by 
\[
\theta_i = \arccos (\mu_i), \quad 1\leq i \leq p.
\]
The distance between $P$ and $Q$ is defined by 
\[
d(P, Q) = \sqrt{\sum \theta_i^{2}}.
\]
It is a natural generalization of the canonical distance of Euclidean sphere.
Thus we have 
\[
\lambda_i = \tan\theta_i.
\]
In the sequel, we shall assume $n\geq m$ without loss of generality. Let $\a = n+\a'$ and denote $\a$ for $\a'$ for simplicity.

Now we fix $P_0 \in G_{n, m}$. We represent it by the $n$-vector $\e_1\wedge\cdot\cdot\cdot \wedge \e_i \wedge \cdot\cdot\cdot \e_n$. We choose $m$ vectors $\e_{n+\a}$, such that $\{ \e_i, \e_{n+\a} \}$ form an orthonormal basis of $\mathbb{R}^{m+n}$. Denote
\[
\mathbb{U}:= \{ P \in G_{n, m}: w(P, P_0) > 0 \}.
\] 
The $v$-function will be 
\[
v(\cdot, P_0) := w^{-1}(\cdot, P_0) \quad {\rm on} \quad \mathbb{U}.
\]
For arbitrary $P \in \mathbb{U}$ determined by an $n\times m$ matrix $Z$, it is easy to see that
\[
v(P, P_0) = \left( \det(I_n +ZZ^T) \right)^{\frac{1}{2}} = \prod_{\a=1}^m \sec \theta_\a = \prod_{\a=1}^m \frac{1}{\mu_\a},
\]
where $\theta_1, \cdot\cdot\cdot, \theta_m$ denotes the Jordan angles between $P$ and $P_0$.

\bigskip

The second fundamental form $B$ of $M^{n}$ in $\mathbb{R}^{m+n}$ is defined by
\[
B_{UW}:= (\overline{\n}_U W)^N
\]
for $U, W \in \Gamma(TM^n)$. We use the notation $( \cdot )^T$ and $(
\cdot )^N$ for the orthogonal  projections into the tangent bundle
$TM^n$ and the normal bundle $NM^n$, respectively. For $\nu \in
\Gamma(NM^n)$ we define the shape operator $A^\nu: TM^n \rightarrow TM^n$
by
\[
A^\nu (U):= - (\overline{\n}_U \nu)^T
\]
Taking the trace of $B$ gives the mean curvature vector $H$ of $M^n$
in $\mathbb{R}^{m+n}$ and
\[
H:= \hbox{trace} (B) = \sum_{i=1}^{n} B_{e_ie_i},
\]
where $\{ e_i \}$ is a local orthonormal frame field of $M^n$.

\vskip24pt

\section{Proof of Theorem \ref{thm1}}

By Corollary 6.2 in \cite{Xin}, the Gauss map of translating solitons is a $V_0^{T}$-harmonic map.
So the study of  translating solitons is naturally related to {\it V-harmonic maps}.
 Recall that a map $u$ from a Riemannian manifold $(M,g)$ to another Riemannian manifold $(N,h)$
is called a $V$-harmonic map if it solves
\[\tau (u)+du(V)=0,\]
where $\tau(u)$ is the tension field of the map $u$, and $V$ is a vector field on $M$ (cf.\cite{CheJosWan15, CheJosQiu12}). Clearly, it is a generalization of the usual harmonic map. Let $\D_V:= \D + \la V, \n \cdot \ra$ and $V=V_0^{T}$.

\vskip12pt

\noindent{\bf Proof of Theorem \ref{thm1}.}
Let $h:= \left( \frac{v}{2-v} \right)^{\frac{3}{2}}$.  Clearly, $h_0= \left( \frac{v_0}{2-v_0} \right)^{\frac{3}{2}}=4$ and $h_1= \left( \frac{v_1}{2-v_1} \right)^{\frac{3}{2}}<4$  are two constants. Choosing a constant $h_2$ such that $h_1 < h_2 <4.$ Since $v\geq 1$, thus we have $1\leq h\leq h_1<h_2<4.$

 Let $\{e_i\}$ be a local orthonormal frame field on $M^n$ such that $\n e_i = 0$ at the considered point.
 From the translating soliton equation (\ref{eqn-T111}), we derive
\begin{equation*}
\n_{e_j} H =  \left( \overline{\n}_{e_j}(V_0-\la V_0, e_k \ra e_k) \right)^N = - \la V_0, e_k \ra B_{e_j e_k}
\end{equation*}
and
\begin{equation*}
\n_{e_i}\n_{e_j} H = - \la V_0, e_k \ra \n_{e_i} B_{e_j e_k} -\la H, B_{e_i e_k} \ra B_{e_j e_k}.
\end{equation*}
Hence using the Codazzi equation, we obtain that
 \begin{equation*}\aligned
\D_V |H|^2 = & \D |H|^2 +\la V, \n |H|^2 \ra  \\
=& 2\la \n_{e_i}\n_{e_i}H, H \ra +2|\n H|^2+ \la V, \n |H|^2\ra \\
=& -2 \la H, B_{e_i e_k} \ra^2 - 2\la \n_{V_0^{T}} H, H \ra +2|\n H|^2 + \la V, \n |H|^2\ra \\
=& -2 \la H, B_{e_i e_k} \ra^2 - \n_V |H|^2 +2|\n H|^2 + \la V, \n |H|^2\ra \\
=&  -2 \la H, B_{e_i e_k} \ra^2  + 2|\n H|^2.
\endaligned
\end{equation*}
It follows that 
 \begin{equation}\label{eqn-MC}\aligned
\D_V |H|^2 \geq   2|\n H|^2-2|B|^2|H|^2.
\endaligned
\end{equation}

For any $X=(x_1, x_2, ..., x_{m+n}) \in \mathbb{R}^{m+n}$, let $r=|X|$,  we have 
\begin{equation}\label{eqn-ER}\aligned
\nabla r^2 =&  2X^{T},  \quad |\nabla r| \leq 1 \\
\Delta r^2 = & 2n + 2\la H,  X \ra \leq 2n+2r.
\endaligned
\end{equation}
Let $B_a(o)$ be the closed ball centered at the origin $o$ with radius $a $ in $\mathbb{R}^{m+n}$ and $D_a(o)=M^n\cap B_a(o)$. Let $\gamma: M^n\to G_{n,m}$ be the Gauss map. Define $f: D_a(o) \rightarrow \mathbb{R}$ by
\[
f=\frac{(a^2-r^2)^2|H|^2}{(h_2- h \circ \gamma)^2}.
\]

Since $\left. f \right|_{\partial D_{a}(o)}=0$, $f$ achieves an absolute maximum in the interior of
$D_{a}(o)$, say $f\leq f(q)$, for some $q$ inside $D_{a}(o)$.  
We may also assume $|H|(q) \neq 0$.
Then from
\begin{equation*}
\label{V3.3add}
          \begin{array}{rcl}
           \ds\vs \nabla f(q)&=&0, \\
       \ds \D_V f(q) & \leq & 0,
        \end{array}
\end{equation*}
we obtain the following at the point $q$:
\begin{equation}\label{3.3}
-\frac{2\nabla r^{2}}{a^{2}- r^{2}} + \frac{\nabla |H|^2}{|H|^2}
 + \frac{2\nabla (h\circ \gamma)}{h_2-h\circ \gamma} = 0,
\end{equation}
\begin{equation}\label{3.4}
- \frac{2 \D_V r^{2}}{a^{2}-r^{2}} - \frac{2 \left| \nabla r^{2}
\right|^{2}}{\left( a^{2} - r^{2} \right)^{2}} + \frac{\D_V
|H|^2}{|H|^2} - \frac{\left|\nabla |H|^2\right|^{2}}{|H|^4} +
\frac{2 \D_V (h\circ \gamma)}{h_2-h\circ \gamma} + \frac{2\left|
\nabla (h\circ \gamma) \right|^{2}}{\left( h_2-h\circ \gamma
\right)^{2}} \leq 0.
\end{equation}
Direct computation gives us
\begin{equation}\label{eqn-GMC}\aligned
|\n |H|^2|^2 = |2\la \n H, H \ra|^2 \leq 4|\n H|^2 |H|^2,
\endaligned
\end{equation}
and 
\begin{equation}\label{eqn-GB}\aligned
|\n (h\circ \gamma)| \leq |d h| |\n \gamma| \leq 3\left( \frac{v_1}{2-v_1} \right)^{\frac{5}{2}}\cdot |B| =: C_1 |B|.
\endaligned
\end{equation}
It follows from (\ref{eqn-MC}) and (\ref{eqn-GMC}) that
\begin{equation}\label{eqn-MB}\aligned
\frac{\D_V|H|^2}{|H|^2} \geq \frac{|\n |H|^2|^2}{2|H|^4} - 2|B|^2.
\endaligned
\end{equation}
From (\ref{3.3}), we obtain 
\begin{equation}\label{eqn-MCRB}\aligned
\frac{|\n |H|^2|^2}{|H|^4} \leq \frac{4|\n r^2|^2}{(a^2-r^2)^2}+ \frac{8|\n r^2||\n(h\circ \gamma)|}{(a^2-r^2)(h_2-h\circ \gamma)} + \frac{4|\n(h\circ \gamma)|^2}{(h_2-h\circ \gamma)^2}
\endaligned
\end{equation}
By (4.6) in \cite{XY} and Corollary 6.2 in \cite{Xin}, we obtain
\begin{equation}\label{eqn-LG}\aligned
\D_V(h\circ \gamma) =& \sum_{i=1}^n {\rm Hess} (h) (d\gamma (e_i), d\gamma(e_i)) + dh(\tau(\gamma) + d\gamma (V)) \\
=& \sum_{i=1}^n {\rm Hess} (h) (d\gamma (e_i), d\gamma(e_i)) \geq 3h|d\gamma|^2 =3h|B|^2.
\endaligned
\end{equation}
Substituting (\ref{eqn-ER}), (\ref{eqn-GB}), (\ref{eqn-MB}), (\ref{eqn-MCRB}), (\ref{eqn-LG}) into (\ref{3.4}), we have 
\begin{equation*}\label{eqn-GJ}\aligned
\left( \frac{3h}{h_2-h\circ \gamma} -1 \right) |B|^2 - \frac{4C_1r}{(a^2-r^2)(h_2-h\circ \gamma)}|B| -\frac{2n+4r}{a^2-r^2} -\frac{8r^2}{(a^2-r^2)^2}  \leq 0
\endaligned
\end{equation*}
Since 
\[
\frac{3h}{h_2-h\circ \gamma} -1 = \frac{3h-(h_2-h\circ \gamma)}{h_2- h\circ \gamma} \geq \frac{4-h_2}{h_2-1}.
\]
Denote $C_2 := \frac{4-h_2}{h_2-1}$. Obviously, $C_2$ is a positive constant. Note that $h\leq h_1<h_2$, the above two inequalities then imply
\begin{equation*}\label{eqn-GJ}\aligned
C_2 |B|^2 - \frac{4C_1r}{(a^2-r^2)(h_2-h_1)}|B| -\frac{2n+4r}{a^2-r^2} -\frac{8r^2}{(a^2-r^2)^2}  \leq 0
\endaligned
\end{equation*}
Note an elementary fact that if $ax^{2}-bx-c \leq 0$ with $a, b, c
$ all positive, then
\begin{equation*}
x \leq \max \{ 2b/a, 2\sqrt{c/a} \}.
\end{equation*}
Therefore, at the
point $q$,
\begin{equation*}\aligned\label{3.16}
|B|^2 & \leq  \max \left\{  \frac{64C^{2}_1 r^{2}}
{C_2^{2}(a^{2}-r^{2})^{2}(h_2-h_1)^{2}},  \frac{4(2n+4r)}{C_2(a^{2}-r^{2})} 
 +  \frac{32r^{2}}{C_2\left( a^{2}-r^{2} \right)^{2}} 
 \right\}.
\endaligned
\end{equation*}
Since $|B|^2 \geq \frac{|H|^2}{n}$, thus we obtain,
 at the point $q$,
\begin{equation}\aligned\label{3.16}
|H|^2  \leq  n\max \left\{  \frac{64C^{2}_1 r^{2}}
{C_2^{2}(a^{2}-r^{2})^{2}(h_2-h_1)^{2}},  \frac{4(2n+4r)}{C_2(a^{2}-r^{2})} 
 +  \frac{32r^{2}}{C_2\left( a^{2}-r^{2} \right)^{2}} 
 \right\}.
\endaligned
\end{equation}
and 
\begin{equation*}\aligned\label{3.16}
f(q)  \leq & n \max \left\{  \frac{64C^{2}_1 a^{2}}
{C_2^{2}(h_2-h_1)^{4}},  \frac{4(2n+4a)a^2}{C_2(h_2-h_1)^2} 
 +  \frac{32a^{2}}{C_2\left( h_2-h_1 \right)^{2}} 
 \right\}.
\endaligned
\end{equation*}
Then for any point $x \in D_{a/2}(o)$, we have
\begin{equation}\aligned\label{3.16'}
|H|^2(x)  \leq & \frac{(h_2-h\circ \gamma)^2}{(a^2-r^2)^2}f(q) \\
\leq & \frac{16n(h_2-1)^2}{9a^4} \max \left\{  \frac{64C^{2}_1 a^{2}}
{C_2^{2}(h_2-h_1)^{4}},  \frac{4(2n+4a)a^2}{C_2(h_2-h_1)^2} 
 +  \frac{32a^{2}}{C_2\left( h_2-h_1 \right)^{2}} 
 \right\}.
\endaligned
\end{equation}
Hence we may fix $x$ and let $a \rightarrow \infty$ in (\ref{3.16'}), we then derive that $H\equiv 0.$ 

If $n\geq 3$, by Theorem 3.1 in \cite{JXY18} (see also Theorem 1.1 in \cite{JXY}), $M^n$ is affine linear.

If $n=2$, let $\{e_1, e_2\}$ be a local tangent orthonormal frame field on $M^2$ and $\{\nu_1, \cdot \cdot \cdot, \nu_m \}
$ a local normal orthonormal frame field on $M^2$ such that $\n e_i$ and $\n \nu_\a = 0$ at the
considered point. Let $V=V^i e_i$ and $H=H^\a \nu_\a$. Then from (\ref{eqn-T111}), we get
\begin{equation}\label{eqn-V0}
V_0=V^T_{0} + V^N_{0} = V + H = V^i e_i + H^\a \nu_\a.
\end{equation}
Thus when $H\equiv 0$, we have $V=V_0$, which is a fixed vector in $\mathbb{R}^{m+2}$ with unit length.

Taking the covariant derivative with respect to $e_j$ in (\ref{eqn-V0}), we obtain
\[
\left( \overline{\n}_{e_j} V^i \right) e_i + V^i h^\a_{ij}\nu_\a +\left( \overline{\n}_{e_j}H^\a  \right) \nu_\a - H^\a h^\a_{ij} e_i = 0,
\]
where $B_{e_i e_j} = h^\a_{ij} \nu_\a$. It follows that
\begin{equation*}\aligned\label{3.16'}
\overline{\n}_{e_j} V^i - H^\a h^\a_{ij} = 0
\endaligned
\end{equation*}
and 
\begin{equation}\aligned\label{eqn-H222}
V^ih^\a_{ij} + \overline{\n}_{e_j}H^\a = 0.
\endaligned
\end{equation}
From (\ref{eqn-H222}) and $H^\a = h^\a_{11} + h^\a_{22}$, we derive
\begin{equation*}\aligned\label{eqn-H333}
h^\a_{11} = &- \frac{1}{|V|^2}\left( V^1 \overline{\n}_{e_1}H^\a - V^2\overline{\n}_{e_2} H^\a - H^\a (V^2)^2 \right), \\
h^\a_{12} = &- \frac{1}{|V|^2}\left( V^1\overline{\n}_{e_2}H^\a +V^2 \overline{\n}_{e_1} H^\a + H^\a V^1V^2 \right), \\
h^\a_{22} = & \frac{1}{|V|^2}\left( V^1\overline{\n}_{e_1} H^\a - V^2 \overline{\n}_{e_2}H^\a + H^\a (V^1)^2 \right). 
\endaligned
\end{equation*}
Then direct computation gives us
\begin{equation}\aligned\label{eqn-H333}
|B|^2 = & \sum_{i, j, \a} (h^\a_{ij})^2 = (h^\a_{11})^2  + 2(h^\a_{12})^2   + (h^\a_{22})^2  \\
=& \frac{1}{|V|^4} \left\{ [V^1 \overline{\n}_{e_1}H^\a - V^2\overline{\n}_{e_2} H^\a - H^\a (V^2)^2]^2 \right. \\
&+ 2 [ V^1\overline{\n}_{e_2}H^\a +V^2 \overline{\n}_{e_1} H^\a + H^\a V^1V^2  ]^2  \\
&\left. + [ V^1\overline{\n}_{e_1} H^\a - V^2 \overline{\n}_{e_2}H^\a + H^\a (V^1)^2 ]^2 \right\} \\
=& \frac{1}{|V|^4} \left( 2|V|^2 |\n H|^2+ 2V^1|V|^2 H^\a\overline{\n}_{e_1}H^\a + |V|^4|H|^2 + 2V^2 |V|^2 H^\a \overline{\n}_{e_2}H^\a \right) \\
=& |H|^2 + \frac{2|\n H|^2}{|V|^2} + \frac{\n_{V}|H|^2}{|V|^2}.
\endaligned
\end{equation}
Since $H\equiv 0$, we have $B\equiv 0$ by (\ref{eqn-H333}).  
Namely, $M^2$ is an affine plane. Hence we complete the proof.
\qed

\end{document}